\documentclass[11pt,a4]{article}
\usepackage{eurosym}
\usepackage{amssymb}
\usepackage{amsfonts}
\usepackage{amsmath}
\usepackage{graphicx}
\usepackage{amsthm}
\usepackage{algorithm}
\usepackage{epic}

    \makeatletter
    \renewcommand*{\@fnsymbol}[1]{\ensuremath{\ifcase#1\or *\or **\or \ddagger\or
       \mathsection\or \mathparagraph\or \|\or **\or \dagger\dagger
       \or \ddagger\ddagger \else\@ctrerr\fi}}
    \makeatother

\newtheorem{rem}{Remark}
\newenvironment{remark}{\begin{rem}\rm}{\end{rem}}
\newtheorem{exmpl}{Example}
\newenvironment{example}{\begin{exmpl}\rm}{\end{exmpl}}
\newtheorem{theorem}{Theorem}

\newtheorem{lemma}{Lemma}

\renewcommand{\blacksquare}{\quad \square}

\title{Forward stable computation of roots of real
      polynomials with only real distinct roots}
\author{Nevena\ Jakov\v{c}evi\'{c} Stor\thanks{Faculty of Electrical
  Engineering, Mechanical Engineering and Naval Architecture, University of
  Split, Rudjera Bo\v{s}kovi\'{c}a 32, 21000 Split,
  Croatia, nevena@fesb.hr}
\ and 
Ivan\ Slapni\v{c}ar\thanks{Faculty of Electrical
  Engineering, Mechanical Engineering and Naval Architecture, University of
  Split, Rudjera Bo\v{s}kovi\'{c}a 32, 21000 Split,
  Croatia, ivan.slapnicar@fesb.hr}}

\begin{document}

\maketitle

\begin{abstract}
As showed in (Fiedler, 1990), any polynomial can be
expressed as a characteristic polynomial of a complex symmetric arrowhead
matrix. This expression is not unique.
If the polynomial is real with only real distinct roots, the matrix can be
chosen real. By using accurate forward stable algorithm for computing
eigenvalues of real symmetric arrowhead matrices from (Jakov\v{c}evi\'{c} Stor, Slapni\v{c}ar, Barlow, 2015),
we derive a forward stable
algorithm for
computation of roots of such polynomials in $O(n^2)$ operations. The algorithm computes each root to
almost full accuracy. In some cases, the algorithm invokes
extended precision routines,
but only in the
non-iterative part. Our examples include numerically difficult problems, like
the well-known Wilkinson's polynomials. Our algorithm compares favourably to
other method for polynomial root-finding, like MPSolve or Newton's method.

\end{abstract}



\section{Introduction and Preliminaries}\label{sec:preliminaries}

Polynomials appear in many areas of scientific computing and engineering.
Developing fast algorithms and reliable implementations of polynomial
solvers are of challenging interest. Famous example by James H. Wilkinson in
1963 \cite{Wil65}, usually referred to as \emph{Wilkinson's polynomial}, is
often used to illustrate difficulties when finding the roots of a
polynomial. The polynomial of order $n$ is defined by a simple formula:
\begin{equation*}  
W_{n}\left( x\right) =\prod\limits_{i=1}^{20}(x-i)=\left( x-1\right) \left(
x-2\right) \cdots \left( x-n\right).
\end{equation*}

For example, the location of the roots of $W_{20}$ is very sensitive to
perturbations in the coefficients, so that in \cite{Wil84}, Wilkinson said:
"Speaking for myself, I regard it as the most traumatic experience in my
career as a numerical analyst." Many methods for finding roots of
polynomials with ever increasing accuracy have been developed since (see for
example \cite{BiFi00}, \cite{Gra08}).

In \cite{Fie90}, Miroslav Fiedler showed that any polynomial can be
expressed as a characteristic polynomial of a complex symmetric arrowhead
matrix. This expression is not unique. If the polynomial is real with only
real distinct roots, the matrix can be chosen real. We have the following
theorem:

\begin{theorem}
\cite[Theorem 3]{Fie90}\label{fiedler} Let $u\left( x\right) $ be a
polynomial of degree $n$,
\begin{equation}  \label{poly}
u\left( x\right) =x^{n}+px^{n-1}+r\left( x\right),
\end{equation}
Let
\begin{equation}  \label{D}
D=\mathrm{diag}(d_{1},\ldots,d_{n-1}),
\end{equation}
where $d_j$ are all distinct and $u\left( d_{j}\right) \neq 0$. Let
\begin{align}
v\left( x\right) & =\prod\limits_{j=1}^{n-1}\left(x-d_{j}\right),  \notag \\
\alpha& =-p-\sum\limits_{j=1}^{n-1}d_{j},  \label{alpha} \\
z&=%
\begin{bmatrix}
\zeta_1 & \zeta_2 & \cdots & \zeta_{n-1}%
\end{bmatrix}%
^T,  \notag
\end{align}
where
\begin{equation}  \label{z}
\zeta_{j}^{2}= \frac{-u(d_{j})}{v^{\prime }( d_{j})} \equiv \frac{-u(d_{j})}{%
\displaystyle\prod_{\substack{ i=1  \\ i\neq j}}^{n-1} ( d_{j}-d_i)}.
\end{equation}
Then the symmetric arrowhead matrix
\begin{equation}  \label{A}
A=\left[
\begin{array}{cc}
D & z \\
z^{T} & \alpha%
\end{array}%
\right],
\end{equation}
has characteristic polynomial $\left( -1\right) ^{n}u\left( x\right)$.%
\newline
\vspace*{3mm}

If $u\left( x\right) $ has only real distinct roots and the $d_{j}$'s
interlace them, then $A$ is real.
\end{theorem}

Fiedler concludes his paper by stating "One can hope to obtain, by some
sophisticated special choice of the numbers $d_j$, stable or even universal
algorithms for solving algebraic equations."\footnote{%
In a report by Corless and Litt \cite{CoLi01}, the matrix $A$ from theorem %
\ref{fiedler} is referred to as \emph{generalized companion matrix not
expressed in monomial basis}. In this case, the basis is the Newton basis.}

In \cite{JSB13} the authors developed a forward stable algorithm for
computing eigendecomposition of a real symmetric irreducible arrowhead
matrix, which is exactly the matrix $A$ given by Theorem \ref{fiedler}.%
\footnote{%
In \cite{JSB13}, the arrowhead matrix $A$ is called ``irreducible'' if $d_j$
are all distinct and $z_j \neq 0$, $j=1,\ldots,n-1$.} More precisely, the
algorithm from \cite{JSB13} computes each eigenvalue and all individual
components of the corresponding eigenvector of a given arrowhead matrix of
floating-point numbers to almost full accuracy in $O(n)$ floating
point-operations, a feature which no other method has.

In this case, we are interested only in the roots of $u$, that is, in the
eigenvalues of $A$ from (\ref{A}), each of which is computed independently
of the others in $O(n)$ operations. This, together with independent
computation of elements of $z$, makes our algorithm suitable for parallel
computing.

In this paper, we propose a new two-step algorithm: given a polynomial $u$
of the form (\ref{poly}) whose coefficients are given floating-point numbers,

\begin{enumerate}
\item compute the generalized companion matrix $A$ from (\ref{A}), where the
elements of $z$ and $\alpha$ need to be computed in double the working
precision, and then

\item compute the roots of $u$ as the eigenvalues of $A$ by using modified
version of the forward stable algorithm $aheig$ from \cite[Algorithm 5]%
{JSB13}.
\end{enumerate}

The organization of the paper is the following. In Section \ref{sec:alg}, we describe
our algorithm named $poly\_aheig$ (POLYnomial roots via ArrowHead
EIGenvalues). In Section \ref{sec:accuracy}, we analyse the accuracy of the
algorithm and give forward error bounds -- in Section \ref{sec:A}, we analyse
the accuracy of the comšputed matrix $A$, and in Section \ref{sec:Ai}, we
analyse the accuracy of the computed inverse of the shifted matrix $A$. In
Section \ref{sec:D}, we discuss possible ways
to find the diagonal elements of the matrix $D$ which interpolate the roots
of $u$. In Section \ref{sec:double}, we discuss different implementations of the double
the working precision, including  extended precision routines from
\cite{Dek71} and the Compensated Horner scheme from \cite[Algorithm
4]{Gra08}. Finally, in Section \ref{sec:examples}, we illustrate our algorithm
with two numerically demanding examples and compare it to
the methods from \cite{BiFi00} and \cite{Gra08}.

\section{The algorithm}

\label{sec:alg}

The eigenvalues of the arrowhead matrix $A$ from (\ref{A}) are the zeros of
the function%
\begin{equation*}
\varphi _{A}\left( \lambda \right) =\alpha -\lambda -z^{T}\left( D-\lambda
I\right) ^{-1}z.  \label{eq:17}
\end{equation*}

The forward stable algorithm for solving EVP of arrowhead matrices \cite%
{JSB13} computes all eigenvalues to almost full accuracy. The algorithm is
based on shift--and--invert strategy. Let $d_{i}$ be the pole which is
nearest to $\lambda $. Let $A_{i}$ be the shifted matrix,
\begin{equation}
A_{i}=A-d_{i}I= \left[
\begin{array}{cccc}
D_{1} & 0 & 0 & z_{1} \\
0 & 0 & 0 & \zeta _{i} \\
0 & 0 & D_{2} & z_{2} \\
z_{1}^{T} & \zeta _{i} & z_{2}^{T} & a%
\end{array}%
\right] \quad  \label{Ai}
\end{equation}%
where
\begin{align}  \label{Dza}
D_{1}& =\mathop{\mathrm{diag}}(d_{1}-d_{i},\ldots ,d_{i-1}-d_{i}),  \notag \\
D_{2}& =\mathop{\mathrm{diag}}(d_{i+1}-d_{i},\ldots ,d_{n-1}-d_{i}),  \notag
\\
z_{1}& = [
\begin{array}{cccc}
\zeta _{1} & \zeta _{2} & \cdots & \zeta _{i-1}%
\end{array}%
] ^{T}, \\
z_{2}& = [
\begin{array}{cccc}
\zeta _{i+1} & \zeta _{i+2} & \cdots & \zeta _{n-1}%
\end{array}%
] ^{T},  \notag \\
a& =\alpha -d_{i}.  \notag
\end{align}
Then,
\begin{equation*}  
\lambda=\frac{1}{\nu}+d_i,
\end{equation*}
where $\nu$ is either largest or smallest (first or last) eigenvalue of the
matrix
\begin{equation}
A_i^{-1}\equiv (A-d_iI)^{-1}= = \left[
\begin{array}{cccc}
D_{1}^{-1} & w_{1} & 0 & 0 \\
w_{1}^{T} & b & w_{2}^{T} & 1/\zeta _{i} \\
0 & w_{2} & D_{2}^{-1} & 0 \\
0 & 1/\zeta _{i} & 0 & 0%
\end{array}%
\right] ,  \label{invAi}
\end{equation}%
where%
\begin{align}
w_{1}& =-D_{1}^{-1}z_{1}\frac{1}{\zeta _{i}},  \notag \\
w_{2}& =-D_{2}^{-1}z_{2}\frac{1}{\zeta _{i}},  \notag \\
b& =\frac{1}{\zeta _{i}^{2}} (
-a+z_{1}^{T}D_{1}^{-1}z_{1}+z_{2}^{T}D_{2}^{-1}z_{2}) .  \label{b_ob}
\end{align}

Notice that all elements of the matrix $A_i^{-1}$ are computed with high
relative accuracy, except that in some cases the element $b$ needs to be
computed in double the working precision (for details see \cite{JSB13}).
Also, the elements of $z$ (the Horner scheme) and $\alpha$ (the trace
preservation formula) of $A$ need to be computed in double the working
precision. Notice that our algorithm requires computation in higher
precision only in the finite part, unlike algorithms from \cite{BiFi00,Gra08}%
, which require usage of higher precision in the iterative part.

The described procedure is implemented in the algorithm $poly\_aheig$.

\begin{algorithm}
\caption{}
$\lambda=\mathbf{poly\_aheig}(u,D)$

\% Computes the roots $\lambda$ of the polynomial $u(x)$ from (\ref{poly}) of
order $n$, with

\% distinct real roots. $D$ is defined by (\ref{D}) and its entries interlace the roots

\% of $u(x)$, see Section \ref{sec:D} for details.
\\

\% Compute the values of $u(x)$ in the interpolating points $d_j$ using double

\% the working precision.

\textbf{for} $j=1:n-1$

\ \ $s_{double}(j)=u(d(j))$

\textbf{end}

\% Compute vector $v$ from Theorem \ref{fiedler}  using double the working precision.

\textbf{for} $j=1:n-1$

\ \ $v_{double}(j) =\prod(d(j)-d(1:j-1,j+1:n-1))$

\textbf{end}

\% compute $\alpha$ from Theorem \ref{fiedler} using double the working precision.

\ \ $\alpha_{double}=-p-\sum\limits_{j=1}^{n-1}d_{j}$

\% compute vector $z$ from Theorem \ref{fiedler} using double the working precision.

\textbf{for} $j=1:n-1$

\ \ $\zeta_{double}(j)=\sqrt{-s_{double}(j)/v_{double}(j)}$

\textbf{end}

\% call modified algorithm $aheig$

\textbf{for} $k=1:n$

\ \ $\lambda(k)=\mathbf{aheig\_mod}(D,z_{double},\alpha_{double},k)$

\textbf{end}

\label{alg1}
\end{algorithm}
\clearpage

\begin{remark}
The algorithm $aheig\_mod$ is a simple modification of the algorithm $aheig$
from \cite[Algorithm 5]{JSB13}. The accuracy of the algorithm $aheig$ is
essentially based on the assumption that all elements of the matrix $%
A_i^{-1} $ from (\ref{invAi}) can be computed with high relative accuracy,
that is, $fl([A_i^{-1}]_{jl})=[A_i^{-1}]_{jl}(1+\kappa_{jl}\varepsilon_M)$,
for some modest $\kappa_{jl}$. For all elements of $A_i^{-1}$ but $b$, this
accuracy is achieved by computing them in standard precision using the
standard precision copies of $z$ and $\alpha$. If, according to the theory
from \cite{JSB13}, the element $b$ needs to be evaluated in double the
working precision, formula (\ref{b_ob}) is evaluated using $z_{double}$ and $%
\alpha_{double}$ in order to obtain full possible accuracy. The details of
the analysis follow.
\end{remark}

\section{Accuracy of the algorithm}\label{sec:accuracy}

The error analysis of the algorithm $aheig$ is given in \cite[Sections 3 and
4]{JSB13}. This analysis assumes that $A$ is the given matrix of
floating-point numbers. Here, however, $A$ is computed by using formulas (%
\ref{poly}-\ref{A}), which must be taken into account. We assume that
computations are performed either in the standard floating-point arithmetic
with the machine precision $\varepsilon_M= 2^{-53}\approx 1.1102\cdot
10^{-16}$ (see \cite[Chapter 2]{Hig02} for details) or with double the
working precision with the machine precision $\varepsilon_M^2=
2^{-106}\approx 1.2326\cdot 10^{-32}$. \footnote{%
Thus, the floating-point numbers have approximately 16 significant decimal
digits. The term ``double the working precision'' means that the
computations are performed with numbers having approximately 32 significant
decimal digits, or with the machine precision equal to $\varepsilon_M^2$.}

Let us first consider the errors in the polynomial evaluation. The standard
method for evaluating polynomial $u(x)$ is the Horner's method \cite[Section
5.1]{Hig02}. Let
\begin{equation}  \label{u}
u(x) =\sum\limits_{i=0}^{n}a_{i}x^{n-i}, \quad a_0\equiv 1,
\end{equation}%
and let%
\begin{equation}\label{condux1}
cond(u,x)=\frac{\sum\limits_{i=0}^{n}\left\vert a_{i}\right\vert \left\vert
x\right\vert ^{n-i}}{\left\vert \sum\limits_{i=0}^{n}a_{i}x^{n-i}\right\vert }=%
\frac{\widetilde{u}(x) }{u(x)}.
\end{equation}%
Notice that $cond(u,x)\geq 1$. Let $Horner(x,u)$ denote the value of $u(x)$
computed in floating point accuracy by the Horner scheme. Then, the relative
error in the computed value is bounded by \cite[Section 5.1]{Hig02}
\footnote{%
In \cite{Hig02,Gra08}, the bounds are expressed in terms of quantities $%
\gamma_k= \frac{k\varepsilon_M}{1-k\varepsilon_M}$. For the sake of
simplicity, we use standard first order approximations $\gamma_k\approx
k\varepsilon_M$.}
\begin{equation*}
\frac{\left\vert u\left( x\right) -Horner\left( u,x\right) \right\vert }{%
\left\vert u\left( x\right) \right\vert }\leq cond(u,x)\times 2n\varepsilon
_{M}.
\end{equation*}

Thus, when $Horner(u,x)$ is evaluated in double the working precision, the
relative error is bounded by
\begin{equation*}
\frac{\left\vert u(x) -Horner_{double}(u,x) \right\vert }{\left\vert u(x)
\right\vert }\leq cond(u,x)\times 2 n\varepsilon _{M}^{2}.
\end{equation*}
Therefore,
\begin{equation}
Horner_{double}(u,x)=(1+\kappa_{x} \varepsilon_{M}^2)u(x),  \label{H1}
\end{equation}
where
\begin{equation}
|\kappa_x| \leq cond(u,x)\times 2 n.  \label{H2}
\end{equation}
Notice that, if $cond(u,x)$ is uniformly bounded,
\begin{equation}  \label{condux}
cond(u,x)\leq \frac{1}{\varepsilon_M},
\end{equation}
then
\begin{equation}  \label{kk}
|\kappa_{x}| \leq 2n.
\end{equation}

Other two options to obtain bounds similar to (\ref{H1},\ref{H2}) are to
evaluate all parts of the respective formulas by using extended precision
routines from \cite{Dek71} or compensated Horner scheme from \cite[Algorithm
4]{Gra08} (see Section \ref{sec:double} for details).

We now consider the accuracy of the computed matrices $A$, $A_{i}$ and $%
A_{i}^{-1}$ from (\ref{A}), (\ref{Ai}) and (\ref{invAi}).

\subsection{Accuracy of $A$} \label{sec:A}

Let $\hat{A}$ denote the matrix $A$ computed according to Algorithm 1,
\begin{equation*}
\hat{A}=\left[
\begin{array}{cc}
D & \hat{z}^{(d)} \\
(\hat{z}^{(d)})^{T} & \hat{\alpha }^{(d)}%
\end{array}%
\right] .  
\end{equation*}%
Here $\hat{z}^{(d)}$ and $\hat{\alpha }^{(d)}$ are computed in
double the working precision which we denote by superscript $(d)$. Let
\begin{equation*}
\hat{z}^{\left( d\right) }=%
\begin{bmatrix}
\hat{\zeta }_{1}^{\left( d\right) } & \hat{\zeta }_{2}^{\left(
d\right) } & \cdots  & \hat{\zeta }_{n-1}^{\left( d\right) }%
\end{bmatrix}%
^{T}.
\end{equation*}%

By combining (\ref{z}) and (\ref{H1}), the standard first order
error analysis in double the working precision, gives
\begin{equation}
\hat{\zeta _{j}}^{\left( d\right) }=\sqrt{\frac{-u(d_{j})(1+\kappa
_{d_{j}}\varepsilon _{M}^{2})}{\displaystyle\prod_{\substack{ i=1 \\ i\neq j
}}^{n-1}(d_{j}-d_{i})(1+\varepsilon _{1})(1+\left( n-3\right) \varepsilon
_{2})}(1+\varepsilon _{3})}(1+\varepsilon _{4}),
\end{equation}%
where $|\varepsilon _{1,2,3,4}|\leq \varepsilon _{M}^{2}$. Therefore,
\begin{equation}
\hat{\zeta _{j}}^{\left( d\right) }={\zeta }_{j}(1+\kappa _{\zeta
_{j}}^{(d)}\varepsilon _{M}^{2}),  \label{zeta}
\end{equation}%
where, by using (\ref{H2}),
\begin{equation*}
\left\vert \kappa _{\zeta _{j}}^{(d)}\right\vert \leq \frac{|\kappa
_{d_{j}}|+\left( n-1\right) }{2}+1
\leq n\cdot cond(u,d_j)+\frac{n+1}{2}
. 
\end{equation*}%
Similarly, applying the standard first order error analysis in double the
working precision to (\ref{alpha}), gives
\begin{equation*}
\hat{\alpha }^{\left( d\right) }={\alpha }(1+\kappa _{\alpha
}^{(d)}\varepsilon _{M}^{2}),  
\end{equation*}%
where
\begin{equation}
\left\vert \kappa _{\alpha }^{(d)}\right\vert \leq \frac{\left\vert
a_1\right\vert +\sum\limits_{j=1}^{n-1}\left\vert d_{j}\right\vert }{%
\left\vert \alpha \right\vert }(n-1)\equiv K_{\alpha }(n-1).  \label{alpha1}
\end{equation}

\subsection{Accuracy of $A_i^{-1}$} \label{sec:Ai}

Let $\hat{A}_i^{-1}$ denote the matrix $A_i^{-1}$ computed according to
Algorithm 1 from the matrix $\hat{A}$. All elements of $A_{i}^{-1}$ but
possibly $b$, are computed in standard precision using the standard
precision copies of $\hat{z}^{(d)}$ and $\hat{\alpha}^{(d)}$. Let $%
\hat{\zeta _{j}}$ and $\hat{\alpha }$ denote $\hat{\zeta _{j}}%
^{\left( d\right) }$ and $\hat{\alpha }^{\left( d\right) }$ rounded to
the nearest standard precision number, respectively. Then
\begin{align}
\hat{\zeta _{j}}& =\zeta _{j}\left( 1+\kappa _{\zeta _{j}}\varepsilon
_{M}\right), \quad j=1,\ldots,n-1, \label{fl_zeta} \\
\hat{\alpha }& =\alpha \left( 1+\kappa _{\alpha }\varepsilon _{M}\right)
,  \label{fl_alpha}
\end{align}
where, by using (\ref{zeta})--(\ref{alpha1}),
\begin{align*}
\left\vert \kappa _{\zeta _{j}}\right\vert &\leq \left(\frac{|\kappa
_{d_{j}}|+\left( n-1\right) }{2}\right)\varepsilon_M+1, \quad j=1,\ldots,n-1,\\
\left\vert \kappa _{\alpha }\right\vert & \leq K_{\alpha
}(n-1)\varepsilon_M+1.
\end{align*}
Further, according to (\ref{H2})-(\ref{kk}), if
\begin{equation}  \label{cudj}
cond(u,d_{j})\leq \frac{1}{\varepsilon _{M}},\quad j=1,\ldots,n-1,
\end{equation}%
then 
(\ref{fl_zeta}) holds with
\begin{equation}  \label{kzj}
\left\vert \kappa _{\zeta _{j}}\right\vert \leq n+2,\quad j=1,\ldots,n-1,
\end{equation}
and if
\begin{equation}  \label{Ka}
K_{\alpha }\leq \frac{1}{\varepsilon _{M}},
\end{equation}%
then (\ref{fl_alpha}) holds with
\begin{equation}  \label{ka}
\left\vert \kappa _{\alpha }\right\vert \leq n.
\end{equation}

For $j\notin \{i,n\}$, similarly as in \cite[Proof of Theorem 4]{JSB13}, the
standard first order error analysis gives
\begin{equation*}
[\hat{A}_{i}^{-1}]_{jj} =fl\left(\frac{1}{d_{j}-d_{i}}\right)=\frac{1}{%
d_{j}-d_{i}}(1+\kappa_{jj}\varepsilon_M), \quad |\kappa_{jj}|\leq 2.
\end{equation*}
Similarly, assuming that (\ref{cudj}) and (\ref{kzj}) hold, for $j\notin
\{i,n\}$ we have
\begin{align*}
[\hat{A}_{i}^{-1}]_{ji}& = fl([\hat{A}_{i}^{-1}]_{ij}) = fl \left(
\frac{-\zeta _{j}(1+\kappa _{\zeta _{j}}\varepsilon _{M})} {(d_{j}-d_{i})
\zeta _{i}(1+\kappa_{\zeta_{i}}\varepsilon _{M})} \right) \\
&=\frac{-\zeta _{j}}{(d_{j}-d_{i})\zeta_{i}}(1+\kappa_{ji}\varepsilon_M),
\quad |\kappa_{ji}|\leq (2n+7).
\end{align*}
Finally,
\begin{align*}
[\hat{A}_{i}^{-1}]_{ni}& =fl([\hat{A}_{i}^{-1}]_{in})=fl\left(\frac{1%
} {\zeta _{i}(1+\kappa _{\zeta _{i}}\varepsilon _{M})}\right) \\
&=\frac{1}{\zeta _{i}}(1+\kappa_{ni}\varepsilon_M), \quad |\kappa_{ni}|\leq
(n+3).
\end{align*}

We now analyze the accuracy of the computed element $b$. Let
\begin{equation}
K_{b}=\frac{|\alpha
|+|d_{i}|+|z_{1}^{T}D_{1}^{-1}z_{1}|+|z_{2}^{T}D_{2}^{-1}z_{2}|}{%
|-a+z_{1}^{T}D_{1}^{-1}z_{1}+z_{2}^{T}D_{2}^{-1}z_{2}|},  \label{K_b}
\end{equation}%
where $D_{1}$, $D_{2}$, $z_{1}$, $z_{2}$ and $a$ are defined by (\ref{Dza}).
We have two cases. First, if
\begin{equation*}
K_{b}\not\gg 1,
\end{equation*}%
then $b$ is computed in standard precision using $\hat{\zeta }_{j}$ and $%
\hat{\alpha }$. Let $\hat{b}$ denote the computed $b$. The standard
first order error analysis of (\ref{b_ob}) gives
\begin{align*}
\hat{b} & =fl\bigg(\frac{1}{\zeta _{i}^{2}(1+\kappa _{\zeta
_{i}}\varepsilon _{M})^{2}}\bigg(\alpha (1+\kappa _{\alpha }\varepsilon
_{M})-d_{i}+\sum\limits_{\substack{ j=1 \\ j\neq i}}^{n-1}\frac{\zeta
_{j}^{2}(1+\kappa _{\zeta _{j}}\varepsilon _{M})^{2}}{d_{j}-d_{i}}\bigg)%
\bigg) \\
& =b(1+\kappa _{b}\varepsilon _{M}),
\end{align*}
where
\[
\quad |\kappa _{b}|\leq
(n+2+\max\{2 \max_{j\neq i }|\kappa _{\zeta _{j}}|,|\kappa _{\alpha }|\})\cdot K_{b}+2|\kappa _{\zeta _{i}}|+3.
\]
Additionally, if (\ref{cudj}) and (\ref{Ka}) hold, then (\ref{kzj}) and
(\ref{ka}) hold, as well, and
\begin{equation*}
|\kappa _{b}|\leq (3n+6)\cdot K_{b}+2n+7.
\end{equation*}

Second, if
\begin{equation*}
K_b\gg 1,
\end{equation*}
then, according to the theory from \cite{JSB13}, the element $b$ needs to be
computed in double the working precision using $\hat{\zeta _{j}}^{\left(
d\right) }$ and $\hat{\alpha } ^{\left( d\right) }$ in order to obtain
full possible accuracy. The standard first order error analysis of (\ref%
{b_ob}) in double the working precision gives
\begin{align*}
\hat{b}^{(d)} &=fl\bigg(\frac{1}{\zeta _{i}^{2}(
1+\kappa_{\zeta_{i}}^{(d)}\varepsilon_{M}^2)^2 } \bigg( \alpha(
1+\kappa_{\alpha}^{(d)}\varepsilon_{M}^2) -d_{i} +\sum\limits_{\substack{ %
j=1  \\ j\neq i}}^{n-1}\frac{\zeta _{j}^{2}(
1+\kappa_{\zeta_{j}}^{(d)}\varepsilon _{M}^2)^2} {d_{j}-d_{i}} \bigg) \bigg)
\\
&=b( 1+\kappa_{b}^{(d)}\varepsilon_{M}^2),
\end{align*}
where
\begin{equation*}
|\kappa_b^{(d)}|\leq (n+2+\max\{2 \max_{j\neq i }|\kappa_{\zeta _{j}}^{(d)}|,|\kappa_{\alpha }^{(d)}|\})\cdot K_{b}+2|\kappa _{\zeta _{i}}^{(d)}|+3.
\end{equation*}

Finally, let
\begin{equation}\label{kappaA}
\kappa_{\hat{A}_i}^{(d)}= \max\{2 \max_{j\neq i }|\kappa_{\zeta _{j}}^{(d)}|,|\kappa_{\alpha }^{(d)}|\}.
\end{equation}
If, in addition to (\ref{cudj}) and (\ref{Ka}),
\begin{equation}  \label{final1}
K_b\leq \frac{1}{\varepsilon_M},
\end{equation}
and
\begin{equation}  \label{final2}
\kappa_{\hat{A}_i}^{(d)} \cdot K_b \leq \frac{1}{\varepsilon_M},
\end{equation}
then
\begin{equation*}
\hat{b}=fl\big(\hat{b}^{(d)}\big)=b( 1+\breve{\kappa}_{b}\varepsilon_{M}),
\end{equation*}
where
\begin{equation*}
|\breve{\kappa}_{b}|\leq n+5.
\end{equation*}

The above results are summarized in the following lemma:

\begin{lemma}\label{lemma1}
Let (\ref{cudj}) and (\ref{Ka}) hold, and let $K_b$ be defined by (\ref{K_b}).
For all non-zero elements of the matrix $A_i^{-1}$ from (\ref{invAi}) computed according to
Algorithm 1 and Remark 1, except for the element $[{A}_i^{-1}]_{ii}$, we have
\[
[\hat{A}_i^{-1}]_{kl}=[{A}_i^{-1}]_{kl}(1+\kappa_{kl}\varepsilon_M), \quad
|\kappa_{kl}|\leq (2n+7).
\]

For the computed element $b=[{A}_i^{-1}]_{ii}$ we have the following: if
$K_b\not \gg 1$, then
\[
\hat{b}= b(1+\kappa_b\varepsilon_M), \quad |\kappa_b|\leq (3n+6)\cdot
K_{b}+2n+7.
\]
If $K_b \gg 1$ and if (\ref{final1}) and (\ref{final2}) hold, then
\[
\hat{b}= b(1+\breve\kappa_b\varepsilon_M), \quad |\breve\kappa_b|\leq n+5.
\]
\end{lemma}

The forward error of the computed roots is bounded as follows:

\begin{theorem} \label{T2}
Let (\ref{cudj}) and (\ref{Ka}) hold, and let $K_b$ be defined by (\ref{K_b}).
Let
\begin{equation*}
\hat{\lambda }=\lambda (1+\kappa _{\lambda }\varepsilon _{M})
\end{equation*}%
be the root of $u(x)$ computed according to
Algorithm 1 and Remark 1. If $K_b \not \gg 1$, then
\begin{equation*}
|\kappa _{\lambda }|\leq 3\sqrt{n}[(3n+6)\cdot
K_{b}+2n+7]+3.18 n\left( \sqrt{n}%
+1\right)+4,
\end{equation*}
and if $K_b \gg 1$ and (\ref{final1}) and (\ref{final2}) hold, then
\begin{equation*}
|\kappa _{\lambda }|\leq (6n+21)\sqrt{n}\,+3.18n\left( \sqrt{n}+1\right)+4.
\end{equation*}
\end{theorem}

\proof
Using the same notation as in \cite[\S 3]{JSB13}, the first summand in the
above bound for $\kappa_\lambda$ follows from \cite[Theorems 5 and 6]{JSB13}%
, while the second summand is the error bound for bisection from \cite[\S 3.1%
]{OlSt90}. $\blacksquare $

\subsection{Choosing $d_{j}$}

\label{sec:D}

Finding values of $d_j$ which interpolate roots is not an easy task.
Articles dealing with computing roots of polynomials usually assume that the
initial approximations of the roots are known (see \cite{Gra08}, \cite{Tis01}%
). Another approach, used in \cite{BiFi00}, is to define the polynomial
neighborhood of $u(x)$ as the set of all polynomials with coefficients
having $d_{in}$ common digits with the corresponding coefficients of $u(x)$,
where $d_{in}$ is predefined input precision. Then, the root neighborhood is
the set of the roots of all polynomials in the polynomial neighborhood of $%
u(x)$.

Since our polynomial is real with only real distinct root our proposal is
simpler. Here are some heuristics:

let $\bar{u}(x)$ be the reverse polynomial of the polynomial $u(x)$ from (%
\ref{u}),
\begin{equation*}
\bar{u}(x)=x^{n}u(1/x)=a_0x^n+a_1 x^{n-1} + a_2 x^{n-2} + \cdots + a_{n-2}
x^{2}+ a_{n-1} x + 1.
\end{equation*}

Since the roots of $\bar u(x)$ are the reciprocals of the roots of $u(x)$,
we have two options for the values $d_j$:

\begin{itemize}
\item[-] use the roots of $u^{\prime}(x)$, or

\item[-] use the reciprocals of the roots of $\bar u^{\prime}(x)$.
\end{itemize}

Depending on the magnitude of the roots, their distribution and relative
gaps, one of the methods, or a combination, is expected to work, see Section %
\ref{sec:examples} for examples.

\subsection{Implementation of double the working precision}

\label{sec:double}

We tested three different implementations of double the working precision:

\begin{itemize}
\item[-] convert all quantities to variable precision by Matlab command
\texttt{sym} with parameter \texttt{'f'}, and then evaluate the respective
formulas
-- this is 300 to 1000 times slower than standard precision.

\item[-] convert all quantities from standard 64 bit \texttt{REAL(8)} to 128
bit \texttt{REAL(16)} in Intel \texttt{ifort} \cite{Int}, and then evaluate the
respective formulas -- this is only 3 times slower,

\item[-] evaluate respective formulas by
using extended precision routines \texttt{add2}, \texttt{sub2}, \texttt{mul2}%
, \texttt{div2}, and \texttt{sqrt2} from \cite{Dek71} -- this is $O(10)$
times slower. In these routines double the working precision is simulated by
keeping each number as a pair consisting of higher and lower part of
mantissa. For example, let
\begin{equation*}
[z,zz]=add2(x,xx,y,yy)
\end{equation*}
where all quantities are floating-point numbers with $t$ binary-digits
mantissa. Then
\begin{equation*}
|z+zz -[(x+xx)+(y+yy)]|\leq (|x+xx|+|y+yy|) 2^{-2(t-1)}.
\end{equation*}
If $xx=0$ and $yy=0$, then (exactly) $z+zz = x+y$. We see that this is nearly
equivalent to using double the working precision (the precision is
$\frac{1}{2}\varepsilon_M^2$
instead of $\varepsilon_M^2$).
\end{itemize}

The evaluation of the polynomial $u(x)$ can also be successfully performed by Compensated
Horner scheme from \cite[Algorithm 4]{Gra08}, where both quantities $h$ and
$c$ from this algorithm must be preserved for subsequent computations by
extended precision routines.

\section{Numerical Examples}

\label{sec:examples}

We illustrate our algorithm with two numerically demanding examples. Here
double the working precision in Algorithm 1 was implemented with extended
precision routines from \cite{Dek71}.

\begin{example}
\label{w18}

The coefficients of Wilkinson polynomial $W_{18}$ are, row-wise,\footnote{%
We use $W_{18}$ since all its coefficients are exactly stored as 64-bit
floating-point numbers.}
\begin{equation*}
\begin{array}{rrr}
1 & -171 & 13566 \\
-662796 & 22323822 & -549789282 \\
10246937272 & -147560703732 & 1661573386473 \\
-14710753408923 & 102417740732658 & -557921681547048 \\
2353125040549984 & -7551527592063024 & 17950712280921504 \\
-30321254007719424 & 34012249593822720 & -22376988058521600 \\
6402373705728000 &  &
\end{array}%
\end{equation*}

In this example the interpolating points $d_j$ can be computed by both ways
described in Section \ref{sec:D}, as roots of $u'(x)$ or as the reciprocals of the
roots of $\bar{u}'(x)$. For example, in the latter case we have
\[
\max{K_b}=214.5 \not \gg 1,\quad  \max_j \{ cond(u,d_j)\}=2.62 \cdot 10^{14},
\quad K_\alpha=26.8,
\]
so by Theorem \ref{T2}, the roots of $W_{18}$ are computed by
Algorithm 1 to (almost) full accuracy, in a forward stable manner.

The roots computed by Matlab \cite{Mat} routine \verb|roots|, MPSolve \cite%
{BiFi00} (with 16 decimal digits), Algorithm \ref{alg1} and Mathematica \cite%
{Wol} with 100 digits of precision (properly rounded to 16 decimal digits),
are, respectively:

{\small
\begin{equation*}
\begin{array}{ccc}
\lambda ^{(roots)} & \lambda ^{(MPSolve)} & \lambda ^{(poly\_aheig, Math)}
\\
18.00001193040660 & 18.00000000000000 & 18 \\
16.99987506992020 & 16.99999999999993 & 17 \\
16.00057853967064 & 15.99999999999455 & 16 \\
14.99841877954789 & 15.00000000000043 & 15 \\
14.00282666587300 & 13.99999999999777 & 14 \\
12.99649084561071 & 12.99999999999819 & 13 \\
12.00308090986650 & 12.00000000000329 & 12 \\
10.99809154207482 & 11.00000000000163 & 11 \\
10.00081885564820 & 9.999999999998594 & 10 \\
8.999776556759201 & 9.000000000000055 & 9 \\
8.000029075840132 & 7.999999999999923 & 8 \\
7.000002735870642 & 7.000000000000000 & 7 \\
5.999998227088450 & 5.999999999999999 & 6 \\
5.000000283698958 & 5.000000000000000 & 5 \\
3.999999981972712 & 4.000000000000000 & 4 \\
3.000000000132610 & 3.000000000000000 & 3 \\
2.000000000018936 & 2.000000000000000 & 2 \\
0.999999999999808 & 1.000000000000000 & 1%
\end{array}%
\end{equation*}%
}

Since for every root, the corresponding quantity $K_b\not \gg 1$, the
algorithm $poly\_aheig$ computes fully accurate roots, using only standard
working precision to compute the corresponding matrix $\hat{A}_i^{-1}$ and its
absolutely largest eigenvalue.

MPSolve
requires input to be defined as integers. Also, MPSolve uses 21 decimal
digits to guarantee and obtain relative accuracy of $10^{-13}$, and it uses
234 decimal digits to guarantee and obtain 30 accurate digits.

The Accurate Newton's method from \cite[Algorithm 6]{Gra08} also computes the
roots of $W_{18}$ to full accuracy as described in \cite[Theorem 6]{Gra08}.
However, the starting points $x_0$ which
satisfy the conditions of \cite[Theorem 6]{Gra08}, must be chosen with
greater care and must be relatively close to the desired root (for example,
$x_0=17.1$ to obtain $\lambda_2=17$, or $x_0=1.1$ to obtain $\lambda_{18}=1$.
Since the Accurate Newton's method takes on average 6 steps to convergence for
each root, it needs approximately $12 n^2$ effective extended precision computations, while
our algorithm needs in this case $5n^2$ extended precision computations
to compute the matrix $\hat A$.

The results for $W_{20}$ are similar.

\end{example}

\begin{example}
\label{IN} Consider the polynomial $u$ of degree $5$ with the coefficients
\begin{equation*}
\begin{array}{r}
1.000000000000000e+00 \\
-2.028240960365167e+31 \\
7.136238463529799e+44 \\
-6.277101735386680e+57 \\
4.181389724724491e+42 \\
-6.189700196426900e+26%
\end{array}%
\end{equation*}
or \texttt{sym(u,'f')}
\begin{equation*}
\begin{array}{r}
1 \\
-20282409603651670423947251286016 \\
713623846352979940529142984724747568191373312 \\
-6277101735386680066937501969125693243111159424202737451008 \\
4181389724724490601097907890741292883247104 \\
-618970019642690000010608640%
\end{array}%
\end{equation*}

In this example the interpolating points $d_j$ are efficiently computed as the reciprocals of the
roots of $\bar{u}'(x)$. The values $d_j$ and $cond(u,d_j)$ from (\ref{condux1})
are given in Table \ref{table1}.

\begin{table}
\center
\begin{tabular}{r|c|c}
$j$ & $d_j$ & $cond(u,d_j)$ \\ \hline
1 & 5.277655813324802e+13 & 4 \\
2&  1.759218604441599e+13 & $3.58\cdot 10^{16}$\\
3&  6.253878705847983e-16 & 12.4\\
4&  2.627905491153268e-16 & 46.4
\end{tabular}
\caption{Interpolating points $d_j$ and $cond(u,d_j)$.}
\label{table1}
\end{table}

For the decreasingly ordered roots of $u$, $\lambda_k$, $k=1,2,3,4,5$, the corresponding quantities $K_b$
from (\ref{K_b}),
$\kappa_{\hat{A}_i}^{(d)}$ from (\ref{kappaA}) and their respective products
from (\ref{final2}), all rounded up, are given in Table \ref{table2}.

\begin{table}
\center
\begin{tabular}{r|c|c|c}
$k$ & $K_b$ & $\kappa_{\hat{A}_i}^{(d)}$ & $\kappa_{\hat{A}_i}^{(d)} \cdot K_b$ \\ \hline
1 & 1 & $3.6\cdot 10^{17}$ & $3.6\cdot 10^{17}$ \\
2 & $3.01\cdot 10^{15}$ & $4.7\cdot 10^{2}$ & $1.42\cdot 10^{18}$ \\
3 & $3.01\cdot 10^{15}$ & $4.7\cdot 10^{2}$ & $1.42\cdot 10^{18}$ \\
4 & $12.6$ & $3.58\cdot 10^{17}$ & $4.48\cdot 10^{18}$ \\
5 & $12.6$ & $3.58\cdot 10^{17}$ & $4.48\cdot 10^{18}$.
\end{tabular}
\caption{Values $K_b$, $\kappa_{\hat{A}_i}^{(d)}$ and $\kappa_{\hat{A}_i}^{(d)} \cdot K_b$.}
\label{table2}
\end{table}

We see that the condition (\ref{final1}) is always fulfilled. Also,
$K_\alpha=1$ from (\ref{alpha1}), so (\ref{Ka}) is fulfilled. The condition
(\ref{final2}) does not hold literally. However, we have
$\kappa_{\hat{A}_i}^{(d)} \cdot K_b\not \gg \frac{1}{\varepsilon_M}$, which is
sufficient to obtain almost full accuracy.

The roots computed by Matlab \cite{Mat} routine \verb|roots|, MPSolve \cite%
{BiFi00} (with 16 decimal digits), Algorithm \ref{alg1} and Mathematica \cite%
{Wol} with 100 digits of precision (properly rounded to 16 decimal digits),
are, respectively: {\small
\begin{equation*}
\begin{array}{ccc}
\lambda ^{(roots)} & \lambda ^{(MPSolve)} & \lambda ^{(poly\_aheig, Math)}
\\
2.028240960365167e+31 & 2.028240960365167e+31 & 2.028240960365167e+31 \\
1.759218604441600e+13 +1.538e+8i & 1.759218604441608e+13 &
1.759218623050247e+13 \\
1.759218604441600e+13 -1.538e+8i & 1.759218604441591e+13 &
1.759218585832953e+13 \\
0 & 4.440892098500624e-16 & 4.440892098500624e-16 \\
0 & 2.220446049250314e-16 & 2.220446049250314e-16%
\end{array}%
\end{equation*}%
}

We see that the roots computed by Algorithm \ref{alg1} coincide fully
with roots computed by Mathematica. Here, in addition to the elements $\hat z^{(d)}$
and $\hat \alpha^{(d)}$ of the matrix $\hat A$, the element $b$ of $\hat
A_2^{-1}$ was computed in double the working precision.

Again, MPSolve requires input to be defined as integers, and it uses 21
decimal digits to guarantee and obtain relative accuracy of $10^{-14}$, and
uses 234 decimal digits to guarantee 30 accurate digits.

Here the  Accurate Newton's method from \cite[Algorithm 6]{Gra08} also computes the
roots to full accuracy, provided the respective starting points are chosen
with greater care. However, the conditions of \cite[Theorem 6]{Gra08} cannot
be used - for example, for the largest root $\lambda_1$, there is no starting
point $x_0$ which satisfies the conditions, except  $\lambda_1$ itself.
For $\lambda_2$, the starting
point $x_0$ which satisfies the conditions can differ from $\lambda_2$ in just
last digit.

\end{example}

\end{document}